
\baselineskip=14pt
\parskip=10pt
\def\Tilde{\char126\relax}
\def\halmos{\hbox{\vrule height0.15cm width0.01cm\vbox{\hrule height
 0.01cm width0.2cm \vskip0.15cm \hrule height 0.01cm width0.2cm}\vrule
 height0.15cm width 0.01cm}}
\font\eightrm=cmr8  
\font\eighttt=cmtt8
\magnification=\magstephalf

\parindent=0pt
\overfullrule=0in
 
\centerline
{\bf There Are MORE THAN 2**(n/17) n-LETTER TERNARY SQUARE-FREE WORDS}
\bigskip
\centerline{ {\it 
Shalosh B. EKHAD{$\,^1$} and Doron ZEILBERGER
}\footnote{$^1$}
{\eightrm  \raggedright
Department of Mathematics, Temple University,
Philadelphia, PA 19122, USA. 
{\eighttt [ekhad, zeilberg]@math.temple.edu} 
{\eighttt http://www.math.temple.edu/\Tilde [ekhad,zeilberg]/ },
Aug. 28, 1998. Supported in part by the NSF.
} 
}
 
{\bf Abstract:} We prove that the `connective constant'
for ternary square-free words is at least $2^{1/17} = 1.0416 \dots $,
improving on Brinkhuis and Brandenburg's lower bounds 
of $2^{1/24}=1.0293 \dots$ and
$2^{1/22}=1.032 \dots$ respectively. This is the first improvement
since 1983.
 
A word is {\it square-free}  if it never stutters, i.e. if it
cannot be written as axxb for words a,b and non-empty word x.
For example, `example' is square-free, but `exampample' is not.
See Steven Finch's famous Mathematical
Constants site[3] for a thorough discussion and
many references.
Let $a(n)$ be the number of ternary square-free $n$-letter words
( {\bf A006156, M2550} 
in the Sloane-Plouffe[4] listing, $1,3,6,12,18,30,42, \dots$).
Brinkhuis[2] and Brandenburg[1] showed that $a(n) \geq 2^{n/24}$,
and $a(n) \geq 2^{n/22}$ respectively.
Here we show, by extending the method of [2], that
$a(n) \geq 2^{n/17}$, and hence that 
$\mu:=\lim_{n \rightarrow \infty} a(n)^{1/n}\geq 2^{1/17}= 1.0416 \dots$.
 
{\bf Definition:} A triple-pair 
$[\,[\,U_0,V_0\,]\,,\,[\,U_1,V_1\,]\,,\,[\,U_2,V_2\,]\,]$
where $U_0, V_0, U_1, V_1, U_2, V_2$ are words in the
alphabet $\{0,1,2\}$ of the same length $k$, will be called
a {\it $k$-Brinkhuis triple-pair} if the following conditions
are satisfied.
 
$\bullet$ The $24$ words of length $2k$,
$$ 
[U\,\,or\,\,V]_0 [U\,\,or\,\,V]_1 \,,
[U\,\,or\,\,V]_0 [U\,\,or\,\,V]_2\,,
[U\,\,or\,\,V]_1 [U\,\,or\,\,V]_2\,,
$$
$$
[U\,\,or\,\,V]_1 [U\,\,or\,\,V]_0\,,
[U\,\,or\,\,V]_2 [U\,\,or\,\,V]_0\,,
[U\,\,or\,\,V]_2 [U\,\,or\,\,V]_1\,,
$$
(i.e.  $U_0U_1$, $U_0V_1, \dots, V_2V_1$), are all square-free.
 
$\bullet$ For every length $r$, $ k/2 \leq r<k$, the $12$ words consisting
of the heads and tails of $\{ U_0,U_1,U_2,V_0,V_1,V_2\}$ of length $r$
are all distinct. \halmos
 
It is easy to see (do it from scratch, or adapt the argument in [2]),
that if 
$[\,[\,U_0,V_0\,]\,,\,[\,U_1,V_1\,]\,,\,[\,U_2,V_2\,]\,]$
is a $k$-Brinkhuis triple-pair,
then for every square-free word
$x=x_1 \dots x_n$  of length $n$ in the alphabet
$\{0,1,2\}$, the $2^n$ words of length $nk$, 
$[U\,\, or\,\, V]_{x_1}[U\,\, or \,\,V]_{x_2} \dots [U\,\, or \,\,V]_{x_n}$ are
also all square-free. Thus the mere existence of
a $k$-Brinkhuis triple-pair implies that
$a(nk) \geq 2^n a(n)$, which implies that $\mu \geq 2^{1/(k-1)}$.
 
{\bf Theorem:} The following is an $18$-Brinkhuis triple-pair
$$
[\,\,[\,\,2 \,\,1 \,\,0\,\, 2 \,\,0 \,\,1 \,\,2 \,\,0\,\, 2 
\,\,1 \,\,2\,\, 0\,\, 1\,\, 0 \,\,2\,\, 0\,\, 1\,\, 2 \quad , \quad
\,\,2\,\,1\,\,0\,\,2\,\,0\,\,1\,\,0\,\,2\,\,1\,\,2\,\,0\,\,2\,\,1
\,\,0\,\,2\,\,0\,\,1\,\,2\,\,],
$$
$$
\,\,\,[\,\,0 \,\,2 \,\,1\,\, 0 \,\,1 \,\,2 \,\,0 \,\,1\,\, 0 
\,\,2 \,\,0\,\, 1\,\, 2\,\, 1 \,\,0\,\, 1\,\, 2\,\, 0 \quad , \quad
\,\,0\,\,2\,\,1\,\,0\,\,1\,\,2\,\,1\,\,0\,\,2\,\,0\,\,1\,\,0\,\,2
\,\,1\,\,0\,\,1\,\,2\,\,0\,\,],
$$
$$
\,\,\,\,\,[\,\,1 \,\,0 \,\,2\,\, 1 \,\,2 \,\,0 \,\,1 \,\,2\,\, 1 
\,\,0 \,\,1\,\, 2\,\, 0\,\, 2 \,\,1\,\, 2\,\, 0\,\, 1 \quad , \quad
\,\,1\,\,0\,\,2\,\,1\,\,2\,\,0\,\,2\,\,1\,\,0\,\,1\,\,2\,\,1\,\,0
\,\,2\,\,1\,\,2\,\,0\,\,1\,\,]\,\,] \, .
$$

{\bf Proof}: Purely Routine! \halmos
 
{\bf Remark:} The above $18$-Brinkhuis triple-pair was found by
the first author by running procedure
{\tt FindPair();} in the Maple package {\tt JAN}, written
by the second author. {\tt JAN} is available from
the second author's website
{\tt http://www.math.temple.edu/\Tilde zeilberg/}
(Click on {\tt Maple programs and packages}, and then on 
{\tt JAN}.)
 
{\bf Another Remark:} Brinkhuis[2] constructed a $25$-Brinkhuis triple-pair
in which $U_0$ and $V_0$ were palindromes, and 
$U_1$, $U_2$, were obtained from $U_0$ by adding, component-wise,
$1$ and $2$ mod $3$, respectively, and similarly for $V_1$, $V_2$. 
Our improved example resulted from relaxing the superfluous condition of
palindromity, but we still have the second property.
It is very likely that by relaxing the second property,
it would be possible to find even shorter Brinkhuis triple-pairs, and
hence get yet better lower bounds for $\mu$. Alas, in this
case the haystack gets much larger!

{\bf References}
 
1. F.-J. Brandenburg, Uniformly growing kth power-free homomorphisms,
       Theor. Comp. Sci. 23 (1983) 69-82.
 
2. J. Brinkhuis, Non-repetitive sequences on three symbols, Quart. J.
       Math. Oxford (2) 34 (1983) 145-149.
 
3. S. Finch, {\it ``Favorite Mathematical Constants Website''},
{
\hfill \break
\tt http://www.mathsoft.com/asolve/constant/words/words.html} .
 
4. N.J.A. Sloane and S. Plouffe, {\it ``The Encyclopedia
of Integer Sequences''}, Academic Press, 1995.
(Online: {\tt http://www.research.att.com/\Tilde njas/sequences/ }),
Direct URL for M2550:
{
\tt
http://www.research.att.com/cgi\Tilde bin/access.cgi/as/njas/sequences/eisA.cgi?Anum=006156.
} 
 
\bye